\documentclass[A4paper,12pt]{article}
\usepackage{amsmath}
\usepackage{latexsym}
\usepackage{mathrsfs}
\usepackage{amssymb}
\usepackage{amscd}
\usepackage[dvips]{graphicx}                  

\newtheorem{theorem}{Theorem}
\newtheorem{corollary}{Corollary}

\newenvironment{proof}[1][Proof]{\textbf{#1.} }{\ \rule{0.5em}{0.5em}}
\date{}
\long\def\symbolfootnote[#1]#2{\begingroup%
	\def\thefootnote{$\;$}\footnote[#1]{$^*$#2}\endgroup}
\begin{document}
	
	\title{A note on Kuratowski partitions \\in metric spaces}
	\author{Joanna Jureczko}
\maketitle

\symbolfootnote[2]{Mathematics Subject Classification: Primary 03C25, 03E35, 03E55, 54E52.

	\hspace{0.2cm}
	Keywords: \textsl{Kuratowski partition, precipitous ideal, K-ideal, real-measurable cardinal, measurable cardinal.}}

\begin{abstract}
We prove that the existence of a complete metric space of cardinality at most  $2^{\kappa}$ admitting Kuratowski partition is a consequence of $\kappa$ being the smallest real-valued measurable cardinal not greater than $ 2^{\aleph_0}$.
\end{abstract}

\maketitle

\section{Introduction}

The notion "Kuratowski partition" has its origin in the old question posed in 1935, (see \cite{KK}), by K Kuratowski, who asked  whether a function $f \colon X \to~Y$, (where $X$ is completely metrizable and $Y$ is metrizable), such that each preimage of an open set of $Y$ has the Baire property, is continuous apart from a meager set.
This problem was later considered on the occassion of investigations around the existence of non-measurable sets, (i.e. \cite{LB} and unpublished results of Solovay). 
In   \cite{EFK} the authors showed that Kuratowski's problem is equivalent to the  existence of partitions of completely metrizable spaces into meager sets such that the union of each subfamily  of this partition has the Baire property. Such a partition is called  \textit{Kuratowski partition}.

In \cite{FK} the authors proved that, if ZFC + "there exists a Kuratowski partition of a Baire metric space" is consistent, then ZFC + "there exists measurable cardinal" is consistent as well, by using forcing methods and the localization property. 

 In \cite{JJ} there was proposed the notion of $K$-ideal which was the previous idea to attact the result solved here. Unfortunately, this idea occured unworkable because based only on the structure of a $K$-ideal one cannot "rebuild" a space admitting Kuratowski partition. 
Moreover, it is necessary to assume that the special space $X(I)$, where $I$ is a $K$-ideal, is complete (see \cite{FK} or the next section for formal definitions) because, as will be shown here, a $K$-ideal $I$ is maximal if $X(I)$ is complete. Thus, the assumption used in the main result of this paper is the only one with which one can show the existence of a complete metric space admitting Kuratowski partition. The restriction of cardinality of a space which is assumed  in the main result cannot be omitted. It follows from \cite{EFK, JJ}.

This paper is the continuation of \cite{JJ2, FJW, JJ}.
The main result of this paper is to prove that the assumption of $\kappa$ being the smallest real-valued measurable cardinal at most $2^{\aleph_0}$ implies the existecce of a complete metric space of cardinality at most $ 2^{\kappa}$ admitting Kuratowski partition. 

The ideas of proofs of Theorem 1 and Theorem 2 are partly based on \cite{FJ1} but in that paper there are incorrect and incomplete. Here, the proofs have been corrected and supplemented. Corollary 2 comes from \cite{FJ1}. 

The results presented here complete the investigations around spaces admitting  Kuratowski partitions
and have wider applications, including some in measurable selector theory, (compare \cite{JJ1}), and related topics.

\section{Definitions and previous results}

This section included definitions and notion used in this paper.  For definitions and facts not cited here we refer interested readers to, e.g., \cite{KK1} (topology) and \cite{TJ} (set theory).
\\
\\
\textbf{2.1.}
Throughout the paper, we assume that $X$ is a Baire space.
A set $U \subseteq X$ has \textit{the Baire property} iff there exists an open set $V \subset X$ and a meager set $M \subset X$ such that $U = V \triangle M$, where $\triangle$ represents the symmetric difference of sets.
\\
\\
\textbf{2.2.}
A partition $\mathcal{F}$ of $X$ into meager subsets of $X$ is called  \textit{Kuratowski partition} iff $\bigcup \mathcal{F}'$ has the Baire property for all $\mathcal{F}' \subseteq \mathcal{F}$. Then we say that 
$X$ \textit{admits Kuratowski partition}.

We define and assume in the whole paper that $$\kappa = \min\{|\mathcal{F}|\colon \mathcal{F} \textrm{ is Kuratowski partition of }X \}.$$
For a given cardinal $\kappa$ we enumerate
$\mathcal{F} = \{F_\alpha \colon \alpha < \kappa\}$. In the whole paper we assume that $\kappa$ is a regular cardinal.

For an open set $U \subseteq X$ treated as a subspace of $X$ that is Baire, the family $$\mathcal{F}\cap U = \{F \cap U \colon F \in \mathcal{F}\}$$ is a Kuratowski partition of $U$.
\\
\\
\textbf{2.3.}
We say that a map $f \colon X \to Y$ has \textit{the Baire property} iff for each open set $V \subset Y$, $f^{-1}(V)$ has the Baire property.
\\

The next fact comes from \cite{FJW} and can be recognized as the characterization of a space admitting Kuratowski partition.
\\\\
\textbf{Fact 1 (\cite{FJW})}  \textit{Let $X, Y$ be topological spaces and $A \subset X$. The following statements are equivalent:}
\begin{itemize}
\item [(i)] \textit{The set $A$ does not admit Kuratowski partition.}
\item [(ii)] \textit{For any mapping $f \colon A \to Y$ having the Baire property, there exists a meager set $M \subset A$ such that $f\upharpoonright(A\setminus M)$ is continuous.}
\end{itemize}

\noindent
\textbf{2.4.}
With any Kuratowski partition
$\mathcal{F} = \{F_\alpha \colon \alpha < \kappa\}$ indexed by  $\kappa$,  we may associate an ideal 
$$I_\mathcal{F} = \{A \subset \kappa \colon \bigcup_{\alpha \in A} F_\alpha \textrm{ is meager}\},$$
which we call  \textit{$K$-ideal}. Obviously, $I_\mathcal{F}$ is non-principal and $[\kappa]^{< \kappa}~\subseteq I_{\mathcal{F}}$.
\\
\\
\textbf{2.5.}
Let $I$ be an ideal on $\kappa$ and let $S$ be a set with positive measure, i.e., $S \in P(\kappa) \setminus I$. (For the convenience, we use $I^+$ instead of $P(\kappa) \setminus I$ throughout). 
\\An \textit{$I$-partition} of $S$ is a maximal family $W$ of subsets of $S$ of positive measure such that $A \cap B \in I$ for all distinct $A, B \in W$.

An $I$-partition $W_1$ of $S$ is a \textit{refinement} of an $I$-partition $W_2$ of $S$ ($W_1 \leq~W_2$) iff each $A \in W_1$ is a subset of some $B\in W_2$.

If $I$ is a $\kappa$-complete ideal on $\kappa$ containing singletons, then $I$ is called \textit{precipitous} iff, whenever $S \in I^+$ and $\{W_n \colon n < \omega\}$ is a sequence of  $I$-partitions of $S$ such that 
$W_0 \geq W_1\geq ... \geq W_n \geq ...$,
there exists a sequence of sets
$X_0 \supseteq X_1\supseteq ... \supseteq X_n \supseteq ...$
such that $X_n \in W_n$ for each $n\in \omega$ and $\bigcap_{n=0}^{\infty} X_n \not = \emptyset$ (see also \cite[p. 424-425]{TJ}).
\\

The next fact shows connections between ideals associated with Kuratowski partitions and precipitousness.
\\\\
\textbf{Fact 2 (\cite{JJ2})} \textit{Let $X$ be a Baire space admitting Kuratowski partition $\mathcal{F}$ of cardinality $\kappa = \min\{|\mathcal{G}| \colon \mathcal{G} \textrm{ is a $K$-partition for  }X\}$.   Then there exists an open set $U \subseteq X$ such that the $K$-ideal $I_{\mathcal{F}\cap U}$ on $\kappa$ associated with $\mathcal{F}\cap U$ is precipitous}.		
\\
\\
\textbf{2.6.}
Let $I$ be an ideal over a cardinal $\kappa$ and let 
$$X(I) = \{x \in (I^+)^\omega \colon \bigcap \{x(n) \colon n \in \omega\} \not = \emptyset \textrm{ and } \forall_{n \in \omega}\ \bigcap\{x(m) \colon m < n\} \not \in I\}.$$
The  set $X(I)$ is considered as a subset of a complete metric space $(I^+)^\omega$, where  $I^+$ is equipped with the discrete topology.
\\

The next fact shows connections between the properties of $X(I)$ and precipitousness of $I$.
\\\\
\textbf{Fact 3  (\cite{FK}) }
\begin{itemize}
	\item [(i)] \textit{$X(I)$ is a Baire space iff $I$ is a precipitous ideal.}
\item [(ii)] \textit{If $I$ be a precipitous ideal over some regular cardinal, then there is Kuratowski partition of $X(I)$.}
\end{itemize}
\noindent
\textbf{2.7.} Let $\lambda$ be a cardinal.
An ideal $I$ is \textit{$\lambda$-saturated} iff there exists no $I$-partition $W$ of size $\lambda$.
Then,
$sat(I)$ \textrm{ is the smallest $\lambda$ such that $I$ is $\lambda$-saturated}.
\\\\
\noindent
\textbf{2.8.}
An uncountable cardinal $\kappa$ is \textit{measurable} iff there exists a non-principal maximal and $\kappa$-complete ideal on $\kappa$.
\\
\\
\textbf{Fact 4 (\cite{TJ})} 
\textit{(i) If $\kappa$ is a regular uncountable cardinal that carries a precipitous ideal, then $\kappa$ is measurable in some transitive model of ZFC}.
\\
\textit{(ii) If $\kappa$ is a  measurable cardinal, then there exists a generic extension in which $\kappa = \aleph_1$, and $\kappa$ carries a precipitous ideal}.
\\
\\
\textbf{Fact 5 (\cite{FK}) } \textit{ZFC +  "there exists measurable cardinal" is equiconsistent with ZFC + "there exists a Baire metric space with a Kuratowski partition of cardinality $\kappa$"}.
\\
\\
\textbf{2.9.} 
A \textit{nontrivial measure} on  $X$ is a map $\mu\colon P(X) \to [0,1]$  such that $\mu$ is a countably additive measure vanishing on points with $\mu(X) =1$ (where $P(X)$ represents the power set of $X$).

Let $\{A_\xi \colon \xi < \kappa\}$ be a family of sets such that $\mu(A_\xi) = 0$ for any $\xi < \kappa$. We say that a measure $\mu$ is \textit{$\kappa$-additive}  $\mu(\bigcup_{\xi < \lambda} A_\xi) = 0$ for any $\lambda < \kappa$.
Let
$$add(\mu) = \min \{\kappa \colon \mu (\bigcup_{\xi < \kappa}A_\xi) > 0, \mu (A_\xi) = 0\}.$$
A cardinal $\kappa$ is \textit{real-valued} iff $\kappa$ carries nontrivial $\kappa$-additive measure.
\\\\
\textbf{Fact 6 (\cite{SU, RS})} \textit{Let $\kappa$ be real-valued measurable. If $\kappa \leqslant 2^{\aleph_0}$, then there is an extension $\mu$ of Lebesgue measure defined on all subsets of $\mathbb{R}$ with $add(\mu) = \kappa$}.
\\\\
\textbf{Fact 7 (\cite{RS})}  \textit{The following theories are equiconsistent}.
\\ \textit{(i) ZFC + "there is a measurable cardinal"}.
\\ \textit{(ii) ZFC + "a Lebesgue measure has a countably additive extension $\mu$ defined on every set of reals"}.
\\\\\textbf{Fact 8 (Ulam, \cite{RS})}  Let $\kappa$ be a real-valued measurable cardinal. Let $\mu$ be a nontrivial  measure on $\kappa$. Then $I = \{A\subseteq \kappa \colon \mu(A) = 0\}$ is $\aleph_1$-saturated.

\section{Main results}

In Fact 2, there is shown that if $X$ is a Baire space  admitting Kuratowski partition $\mathcal{F}$, then there exists an open set $U \subseteq X$ such that  the $K$-ideal $I_{\mathcal{F}\cap U}$ is precipitous. In Theorem 1 we show that  $I_{\mathcal{F}\cap U}$ can be maximal if $X(I_\mathcal{F})$ is complete. Note that this assumption cannot be omitted.

\begin{theorem} If $X$ is a Baire space with Kuratowski partition $\mathcal{F}$ of cardinality $\kappa$, $I_{\mathcal{F}}$ is K-ideal associated with $\mathcal{F}$ and $X(I_{\mathcal{F}})$ is complete, then  $I_{\mathcal{F}\cap U}$ is maximal for some open set $U \subseteq X$.
	\end{theorem}

\begin{proof}
	Let $\mathcal{F} = \{F_\alpha \colon \alpha < \kappa\}$ be Kuratowski partition of $X$ and let $I_\mathcal{F}$ be $K$-partition associated with $\mathcal{F}$. Consider $X(I_\mathcal{F})$. Since $X(I_\mathcal{F})$ is assumed to be complete hence the Baire Theorem holds.
	By Fact 3, $I_\mathcal{F}$ is precipitous. Without loss of generality we can assume that $I_\mathcal{F}$ is everywhere precipitous, i.e. $I_{\mathcal{F}\cap U}$ is precipitous for any open set $U \subset X$.
	
Let $\mathcal{U}$ be a  family of open and pairwise disjoint subsets of $X$, such that $\bigcup \mathcal{U}$ is dense in $X$. Suppose in contrary that $I_{\mathcal{F}\cap U}$ is not maximal for any open  $U\in \mathcal{U}$. 
Without loss of generality we can assume that $sat(I_{\mathcal{F}\cap U})> \omega$ for any $U \in \mathcal{U}$, (compare Fact 8).

For any $U \in \mathcal{U}$ consider a set $S^U \not \in I_{\mathcal{F}\cap U}$ in such a way that $S^U\cap S^{U'} = \emptyset$, whenever $U, U' \in \mathcal{U}$ are different.
For each $ U \in \mathcal{U}$ choose  a sequence of  $I_{\mathcal{F}\cap U}$-partitions of $S^U$  
$$W^U_0 \geq W^U_1\geq ... \geq W^U_n \geq ...$$ and a decreasing sequence 
$$A^U_0 \supset A^U_1\supset ... \supset A^U_n \supset ...$$
with $A^U_k \in W^U_k, (k \in \omega)$
and
$\bigcap_{k=m}^{n} A^U_k \not \in I_{\mathcal{F}\cap U}, (m < n),$ but $\bigcap_{k=0}^{\infty} A^U_k \in I_{\mathcal{F}\cap U}$. (Such a choice is always possible).
Now, define a decreasing sequence  $$A_k = \bigcup_{U \in \mathcal{U}} A_k^U.$$ Obviously, $\bigcap_{k=m}^{n}A_k \not \in I_{\mathcal{F}}$ for any $m<n$ but $\bigcap_{k=0}^{\infty} A_k \in I_\mathcal{F}$. 

Define a sequence $(x_m)_{m \in \omega}$ of elements of $X(I_{\mathcal{F}})$ such that  $x_m = (x_m(k))_{k \in \omega}$ and
$$
x_m(k) = \left\{ \begin{array}{ll}
A_k & \textrm{for } k \leqslant m \\
C   & \textrm{for } k>m\\
\end{array} \right.,
$$
where $C\not \in I_\mathcal{F\cap U}$ is chosen in this way that the following condition holds
 $$d(x_m, x_n) \leqslant \frac{1}{2^{n-m}},$$ for all $m< n$ and $ m, n \in \omega$, ($d$ means metric in $X(I_{\mathcal{F}})$). 

Clearly, the sequence $(x_m)_{m \in \omega}$ fulfills the Cauchy condition, but the limit of  $(x_m)_{m \in \omega}$ does not belong to $X(I_{\mathcal{F}})$, because $\bigcap_{k =0}^{\infty} A_k \in I_{\mathcal{F}}$. A contradiction with the assumption that $X(I_\mathcal{F})$ is complete.
Thus, there is $U \in \mathcal{U}$ such that $I_{\mathcal{F}\cap U}$ is maximal.
\end{proof}
\\

Applying Theorem 1 and Fact 2, we immediately obtain the following result.

\begin{corollary} Let $X$ be a Baire space admitting Kuratowski partition $\mathcal{F}$ of cardinality $\kappa$, where $\kappa = \min\{|\mathcal{G}|\colon \mathcal{G} \textrm{ is Kuratowski partition of }X \}$. If  $X(I_{\mathcal{F}})$ is complete, then  $\kappa$ is measurable.
\end{corollary}

In Theorem 2 we show that there is a complete metric space admitting Kuratowski partition. We emphasize that Theorem 2 is true in ZFC only, and the given space has to have cardinality not greater than $2^{\aleph_0}$ (see \cite{EFK}).

\begin{theorem}
If $\kappa$ is a regular  and  the smallest real-valued measurable cardinal such that $\aleph_1 < \kappa \leqslant 2^{\aleph_0}$, then there exists a complete metric  space of cardinality not greater than $2^{\kappa}$ which admits Kuratowski partition.
\end{theorem}

\begin{proof} For simplification the notation we can assume that $X=[0,1]$. \\
	Let $\mu\colon P(X) \to X$ be nontrivial $\kappa$-additive measure. By Fact 6,  $\mu$ extends everywhere Lebesgue measure  on $X$. 
	Take $\mu$-measurable sets  $A, B \in P(X)$ and define the equivalence relation 
	$$A \sim B \textrm{ iff } \mu(A\triangle B) = 0,$$
	where $\triangle$ means the symmetric difference of sets. 
Let $[A]$  denotes the equivalence class determined by  $A$.

Let
$$Y = \{[A] \colon A \in P(X), A \textrm{ is } \mu\textrm{-measurable}\}.$$
	Define a metric
	$\rho([A], [B]) = \mu(A\triangle B)$ on $Y$.
	Since $A, B \in P(X)$ are $\mu$-measurable, $\rho$ is well defined.	
	Moreover, the space $(Y, \rho)$ is complete because the limit of any sequence $([A_n])_{n \in \omega}$ fulfilling the Cauchy condition is of the form  $[\bigcap_{n \in \omega}\bigcup_{k=0}^{n}A_{k}]$ and belongs to $Y$.
	
	Enumerate the elements of $X$ by $\{x_\alpha \colon \alpha < \mathfrak{c}\}$. For any $x_\alpha \in X$ choose its neighbourhood $U_{x_\alpha}$ in such a way that $U_{x_\alpha} \cap U_{x_\beta} = \emptyset$, whenever $\alpha \not = \beta$.
	 Now, for any $\alpha < \mathfrak{c}$, define 
	$$G_\alpha =\{[A]\in Y \colon \alpha = \min \{\beta < \mathfrak{c} \colon U_{x_\beta} \cap A \not = \emptyset\}\}.$$
	Obviously $G_\alpha \cap G_\beta = \emptyset$ for any pairwise different $\alpha, \beta < \mathfrak{c}$ and $G_\alpha$ is meager in $Y$, (because $Y$ as a complete metric space is a Baire space). 
	To complete the proof it is enough to show that $\bigcup\{G_\alpha \colon \alpha \in D\}$ has the Baire property for any $D \subseteq \mathfrak{c}$. To do this it is enough to show that $\bigcup\{G_\alpha \colon \alpha \in D\}$ contains an open and dense set.

	Let $\mathcal{G} = \{G_\alpha \colon \alpha < \mathfrak{c}\}$ and let $$I_\mathcal{G} = \{D \subseteq \mathfrak{c} \colon \bigcup_{\alpha \in D} G_\alpha \textrm{ is meager}\}.$$  
	Take $D \subseteq \mathfrak{c}$ and consider two cases.
	
	Case 1. If $D \not \in I_{\mathcal{G}}$, define	
	 $$W(D)=\{[A] \in Y \colon \exists_{\alpha \in D}\  U_{x_\alpha} \cap A \not = \emptyset\}.$$
	 Observe that $W(D)$ is contained in $\bigcup\{G_\alpha \colon \alpha \in D\}$ and is open and dense. To see this, take arbitrary $[A] \in W(D)$.
	 Then  $U_{x_\alpha} \cap A \not = \emptyset$ for some $\alpha \in D$ and $\min \{\alpha \in D \colon U_{x_\alpha} \cap A \not = \emptyset\} \in D$. Hence $[A] \in \bigcup\{G_\alpha \colon \alpha \in D\}$.
	 
	 Case 2. If $D \not \in I_{\mathcal{G}}$,
	 then  there exists an open and dense $G_\delta$ - set $$W(D')=\{[A]\in Y \colon \exists_{\alpha \not \in D}\ U_{x_\alpha} \cap A\not = \emptyset\}$$ which is contained in $\bigcup\{G_\alpha \colon \alpha \not \in D\}$. It means that $\bigcup\{G_\alpha \colon \alpha \in D\}$ has the Baire property for any $D \subseteq \mathfrak{c}$.
	Finally, $\mathcal{G}$ is a required Kuratowski  partition of $Y$.
\end{proof}
\\

As the consequence of Fact 1, Fact 5 and Theorem 2 we obtain the following corollary. 

\begin{corollary}
The following theories are consistent:
\\
(1) ZFC + "there is a measurable cardinal",
\\
(2) ZFC + "there is a complete metric space $X$ of cardinality not greater than $2^\mathfrak{c}$ and a function $f \colon X \to Y$ having the Baire property such that there is no meager set $M \subseteq X$ for which $f\upharpoonright(X \setminus M)$ is continuous".
\end{corollary}

{\sc Joanna Jureczko}
\\
Wroc\l{}aw University of Science and Technology, Wroc\l{}aw, Poland
\\
{\sl e-mail: joanna.jureczko@pwr.edu.pl}

\end{document}